\documentstyle[12pt]{article}
\begin{document}

\title{ \Large {  Almost Complex Structure on  $S^{2n}$}
 \footnotetext{ \hskip -0.6cm The author was partially supported by NNSF(10871142) of
 China. \\
 {\it Keywords}. \ Complex structure, Riemannian connection,  Chern class. \\
 {\it Subject classification}. \  53B35, 53C05, 57R20.\\
  E-mail: \ jwzhou@suda.edu.cn }}
\author{Jianwei Zhou}
\date{\small Department of Mathematics, Suzhou University, Suzhou
215006, P.R. China }

\maketitle
\begin{abstract}  We show that there is no  complex
structure in a neighborhood of the space of orthogonal almost
complex structures on the sphere $S^{2n}, \ n>1$.  The method is  to
study the first Chern class of vetcor bundle $T^{(1,0)}S^{2n}$.

\end{abstract}

\baselineskip 15pt
\parskip 3pt

In  [4], we   showed that  the orthogonal twistor space
$\widetilde{\cal J}(S^{2n})$ of the sphere $S^{2n}$ is a Kaehler
manifold and  an orthogonal almost complex structure $J_f$ on
$S^{2n}$ is integrable if and only if the corresponding section
$f\colon\; S^{2n}\to \widetilde{\cal J}(S^{2n}) $ is holomorphic.
These shows there is no  orthogonal complex structure on the sphere
$S^{2n}$ for $n>1$.

We can show that the first  Chern class of the   vector bundle
$T^{H(1,0)} \widetilde{\cal J}(S^{2n})$ can be represented by the
Kaehler form of $ \widetilde{\cal J}(S^{2n})$  with a constant as
coeifficient. Thus $c_1(T^{(1,0)}S^{2n})=f^*c_1(T^{H(1,0)}
\widetilde{\cal J}(S^{2n}))\in H^2(S^{2n})$ is non-zero if the map
$f\colon\; S^{2n}\to \widetilde{\cal J}(S^{2n}) $ is holomorphic. In
following we show that there is no complex structure in a
neighborhood  of the twistor space $\widetilde{\cal J}(S^{2n})$, the
method is  to study $c_1(T^{(1,0)}S^{2n})$.

 Let $\langle \, , \, \rangle$
be the canonical Riemanian metric on the sphere $S^{2n}$ and $J_f$
be an almost complex structure.
$$\mbox{d}s^2_f(X, Y)= \langle X, Y \rangle_f= \frac 12 \langle X, Y \rangle+
\frac 12 \langle J_fX, J_fY \rangle, \ \ \ X,Y\in TS^{2n}\otimes
{\bf C}$$ defines a Hermitian metric on $TS^{2n}\otimes {\bf C}$.
Let  $\tilde e_1,\cdots, \tilde e_{2n}$ be  local $\langle \, , \,
\rangle$-orthonormal frame fields on $S^{2n}$, $\tilde \omega^1,
\cdots, \tilde \omega^{2n}$ be their dual. The almost complex
structure $J_f$ can be represented by $J_f=\sum \ \tilde
e_iJ_{ij}\tilde \omega^j, \ J_f(X)=\sum \ \tilde e_iJ_{ij}\tilde
\omega^j(X)$, the Hermitian metric can be represented by
$$\mbox{d}s^2_f=
\frac 12\sum \ (\delta_{ij} +\sum \ J_{ki}J_{kj})\tilde
\omega^i\otimes \tilde \omega^j.$$

{\bf Remark} \ Let $(P_{ij})$ be the positive definite symmetric
matrix such that $(P_{ij})^2=\frac 12(\delta_{ij} +\sum \
J_{ki}J_{kj})$. The tensor field $P=\sum \ \tilde e_iP_{ij}\tilde
\omega^j$  is well defined. For  any $\langle \, , \,
\rangle_f$-orthogonal almost complex structure $J_1$ and $X,Y\in
TS^{2n}$, we have $$\langle PX, PY\rangle=
  \langle X, Y \rangle_f,$$
 $$\langle PJ_1P^{-1}X, PJ_1P^{-1}Y\rangle =\langle J_1P^{-1}X, J_1P^{-1}Y\rangle_f
=  \langle X, Y \rangle.$$ These shows the almost complex structure
$PJ_1P^{-1}$   preserves  the metric $\langle \, , \, \rangle$.

By the proof of Theorem 2.4 in [4], there are local $\langle \, , \,
\rangle$-orthonormal frame fields $\bar e_1, \bar e_2,\cdots, \bar
e_{2n}$ and their dual $\bar \omega^1,\cdots,\bar \omega^{2n}$ such
that
$$J_f = \sum \ \sqrt{1+\lambda_i^2}
 (\bar e_{2i}\bar \omega^{2i-1} - \bar e_{2i-1}\bar \omega^{2i}) +
  \sum \ \lambda_i (\bar e_{2i-1}\bar \omega^{2i-1} -
\bar e_{2i}\bar \omega^{2i}).$$ Set $\tilde e_{2i-1}=\frac {\sqrt
2}2(\bar e_{2i-1}+\bar e_{2i}), \ \tilde e_{2i}=\frac {\sqrt
2}2(-\bar e_{2i-1}+\bar e_{2i})$, we have
$$J_f =  \sum \ (\lambda_i+\sqrt{1+\lambda_i^2})
 \tilde e_{2i}\tilde \omega^{2i-1}+ \sum \ (\lambda_i-\sqrt{1+\lambda_i^2})
 \tilde e_{2i-1}\tilde \omega^{2i} .$$
Then we have $P=\sum\limits_{k=1}^{2n} \ \mu_{k}\tilde e_k\tilde
\omega^k,$ where
$$\mu_{2i-1}=\left(1+\lambda_i^2+\lambda_i\sqrt{1+\lambda_i^2}\right)^{\frac
12}, \
\mu_{2i}=\left(1+\lambda_i^2-\lambda_i\sqrt{1+\lambda_i^2}\right)^{\frac
12}.$$ Set $e_k=\frac 1 {\mu_k}\tilde e_k$, we have $J_f (e_{2i-1}
)=e_{2i}, \ (PJ_fP^{-1})(\tilde e_{2i-1})= \tilde e_{2i}. \ \ \ \ \
\ \Box$

 Let $\nabla^f$ be the Riemannian connection of the metric
$\langle \, , \, \rangle_f$ and $  e_1,\cdots, e_{2n}$ be
$J_f$-frame fields, $\langle    e_{i} , e_{j}\rangle_f=\delta_{ij},
\ J_f e_{2i-1}= e_{2i},
 \   \omega^1,\cdots,   \omega^{2n}$ be their dual,
$J_f  \omega^{2i-1}=-   \omega^{2i}$. Let
$$\nabla^f (  e_1,\cdots,   e_{2n-1},  e_2,\cdots,   e_{2n})^t=  \omega
(  e_1,\cdots,   e_{2n-1},  e_2,\cdots,   e_{2n})^t,$$ where
$\omega=\left( \begin{array}{cccccccc}
A  & B  \\
C & D
\end{array}\right)$ be the
connection matrix,  $ \Omega=\mbox{d}\omega - \omega\wedge \omega$
 the curvature matrix.   $B+C=A-D=0$ if $J_f$  is integrable and $\nabla^f$ a Kaehler
 connection.

{\bf  Lemma 1} \ $J_f$ is a complex structure if and only if
$J_f(B+C)=A-D$.

{\bf Proof} \ Let $Z_i=  e_{2i-1}-\sqrt{-1}  e_{2i}, \ Z_{\bar i}
=e_{2i-1}+\sqrt{-1}  e_{2i}$ be the $(1,0)$ and $(0,1)$ frame fields
on $S^{2n}$, $i=1,\cdots,n$. Then we have
\begin{eqnarray*}\nabla^f
\left( \begin{array}{cccccccc} Z_i \\ Z_{\bar i}\end{array}\right)
&=& \nabla^f\left( \begin{array}{cccccccc}
I  & -\sqrt{-1}I  \\
I & \sqrt{-1}I
\end{array}\right)\left( \begin{array}{cccccccc}
e_{2i-1} \\ e_{2i}\end{array}\right)\\
&=& \frac 12\left( \begin{array}{cccccccc}
I  & -\sqrt{-1}I  \\
I & \sqrt{-1}I
\end{array}\right)  \left( \begin{array}{cccccccc}
A  & B  \\
C & D
\end{array}\right)\left( \begin{array}{cccccccc}
I  &  I \\ \sqrt{-1}I & -\sqrt{-1}I
\end{array}\right) \left( \begin{array}{cccccccc}
Z_i \\ Z_{\bar i}\end{array}\right)\\
&=&\frac 12\left( \begin{array}{cccccccc}
A+D+\sqrt {-1}(B-C)  & A-D-\sqrt{-1}(B+C) \\
A-D+\sqrt {-1}(B+C)  & A+D-\sqrt{-1}(B-C)
\end{array}\right)\left( \begin{array}{cccccccc}
Z_i \\ Z_{\bar i}\end{array}\right).\end{eqnarray*}

The vector bundle $T^{(1,0)}S^{2n}$ is  generated by
$Z_1,\cdots,Z_n$. As shown in [1], [4], the almost complex structure
$J_f$ is integrable if and only if $\nabla^f_{X}Y\in \Gamma
(T^{(1,0)}S^{2n})$ for any $X, Y \in \Gamma (T^{(1,0)}S^{2n})$. Then
$J_f$ is integrable if and only if $A-D-\sqrt{-1}(B+C)$ are formed
by $(0,1)$-forms, that is
$$J_f(B+C)=A-D. \ \ \  \ \ \ \Box$$

{\bf  Lemma 2} \  The first Chern class of $T^{(1,0)}S^{2n}$ can be
represented by
$$c_1(T^{(1,0)}S^{2n})
=-\frac {1}{4\pi}\mbox{tr} ( \Omega J_0+  \omega \wedge\omega
J_0),$$ where $J_0=\left(
\begin{array}{cccccccc}
{}  & -I  \\
I & {}
\end{array}\right)$.

{\bf Proof} \ With notations used above, the induced connection
$\nabla^f$ on $T^{(1,0)}S^{2n}$ is
$$\nabla^fZ_i
= \sum \ \psi_i^k Z_k,$$ where $\psi_i^k=\frac {1}{2}\sum \ (
\omega_{2i-1}^{2k-1}+ \omega_{2i}^{2k} +\sqrt{-1}(
\omega_{2i-1}^{2k}- \omega_{2i}^{2k-1}))$ and the curvature forms
$\Psi_i^k=\mbox{d} \psi_i^k-\sum \ \psi_i^j\wedge\psi_j^k$,
$\psi_i^k+\overline\psi^i_k =0, \ \Psi_i^k+\overline\Psi^i_k =0$.
 The first Chern class of the vector bundle
$T^{(1,0)}S^{2n}$ can be represented by
$$c_1(T^{(1,0)}S^{2n})=\frac {\sqrt{-1}}{2\pi}\sum \ \Psi_i^i.$$
By $\sum \ \psi_i^j\wedge\psi_j^i=-\sum \ \psi_j^i\wedge\psi^j_i=0$
and $ \sum \ \mbox{d} \omega_{2i-1}^{2i}= \frac 12\mbox{tr}(\mbox{d}
\omega J_0)=\frac {1}{2}\mbox{tr} ( \Omega J_0+  \omega \wedge\omega
J_0)$, we have
$$c_1(T^{(1,0)}S^{2n})=
-\frac {1}{4\pi}\mbox{tr} ( \Omega J_0+  \omega \wedge\omega J_0). \
\ \  \ \ \ \Box$$

As [4], let ${\cal J}(S^{2n})$ be the twistor space on $S^{2n}$, its
sections are the  almost complex structures on $S^{2n}$.

{\bf  Theorem 3} \ When $n>1$, there is no complex structure in a
neighborhood of the space $ \widetilde{\cal J}(S^{2n})$.

{\bf Proof} \ By $$\omega +J_0 \omega J_0=\left(
\begin{array}{cccccccc}
A-D  & B+C  \\
B+C & -A+D
\end{array}\right),\ J_0 \omega - \omega J_0=\left( \begin{array}{cccccccc}
-B-C  & A-D  \\
A-D & B+C
\end{array}\right),$$ we see that the equation
$J_f(B+C)=A-D$ is equivalent to  $J_f( \omega +J_0  \omega
J_0)=J_0(\omega +J_0 \omega J_0).$ Then, if $J_f$ is integrable, we
have
\begin{eqnarray*} && \mbox{tr}( \omega\wedge \omega J_0) \\
&&=\frac 14\mbox{tr}[( \omega +J_0 \omega J_0)\wedge( \omega J_0-J_0 \omega)] \\
&&= -\mbox{tr}[(B+C)\wedge J_f(B+C)] \\
&&= \mbox{tr}[(B+C)\wedge J_f(B+C)^t].
\end{eqnarray*}

The sectional curvature of the metric $\langle \, , \, \rangle$ on
$S^{2n}$ is constant,  $ \Omega_k^l=- \omega^k\wedge \omega^l.$ Then
if $J_f$ is a  $\langle \, , \, \rangle$-orthogonal complex
structure, we have
\begin{eqnarray*}
&& \mbox{tr}( \Omega J_0)+\mbox{tr}( \omega\wedge \omega J_0)   \\
 &  &=\sum \ 2 \omega^{2i-1}\wedge J_f \omega^{2i-1}
 + \sum \ ( \omega_{2i-1}^{2j}+ \omega_{2i}^{2j-1})\wedge J_f
( \omega_{2i-1}^{2j}+ \omega_{2i}^{2j-1}).
\end{eqnarray*}
For any $X \in TS^{2n}$, we have
\begin{eqnarray*}
&& [\mbox{tr}( \Omega J_0)+\mbox{tr}( \omega\wedge \omega J_0)](X, J_fX)   \\
 &  &=-\sum \ 2 ([\omega^{2i-1}(X)]^2+[\omega^{2i}(X)]^2) \\
 && \quad
- \sum \ [( \omega_{2i-1}^{2j}+ \omega_{2i}^{2j-1})(X)]^2 - \sum \
[( \omega_{2i-1}^{2j-1}- \omega_{2i}^{2j})(X)]^2.
\end{eqnarray*}
Then  2-form $\mbox{tr}(\Omega J_0+ \omega\wedge \omega J_0)$ are
non-degenerate everywhere and $S^{2n}$ becomes a symplectic
manifold, this contradict to the fact of $H^2(S^{2n})=0$ for $n>1$.
As the Riemannian curvature is continuous with  the Riemannian
metric, these shows there is a neighborhood of $\widetilde{\cal
J}(S^{2n})$ in ${\cal J}(S^{2n})$ such that there is no complex
structure in this neighborhood. \ \ \ \ \ $\Box$

\vskip 1cm \noindent{\large \bf References} \vskip 0.3cm

{\small

\noindent [1] \  Albuquerque, R.,  Salavessa, I. M. C.: On the
twistor space of pseudo-spheres. Differ. Geom. Appl., {\bf
25}(2007), 207-221.

\noindent [2] \  O'Brian, N. R.,  Rawnsley, J.: Twistor spaces. Ann.
of Global Analysis and Geometry. {\bf 3}(1985), 29-58.

\noindent [3] \  Zhou, J. W.: Grassmann manifold $G(2,8)$ and
complex structure on $S^{6}$.  arXiv: math. DG/0608052.

\noindent [4] \  Zhou, J. W.: The complex structure on  $S^{2n}$.
arXiv: math. DG/0608368.

 \noindent [5] \  Zhou, J. W.: Lectures on differential geometry (in Chinese).
 Science Press, Beijing, 2010.}

\end{document}